\documentclass[12pt,a4paper,oneside,english]{amsart}
\usepackage{babel}
\usepackage{amsmath}
\usepackage{setspace}
\doublespace
 \makeatletter
 \theoremstyle{plain}
 \newtheorem{thm}{Theorem}[section]
 \numberwithin{equation}{section} 
 \numberwithin{figure}{section} 
  \newtheorem*{thm*}{Theorem}

 \newtheorem*{defn*}{Definition}
 \newtheorem{prop}[thm]{Proposition} 
 \newtheorem{rem}[thm]{Remark}
 \newtheorem{claim}[thm]{Claim}
 \newtheorem{lem}[thm]{Lemma} 

\newtheorem*{acknowledgment}{Acknowledgment}

\begin{document}

\title{Rewriting systems in alternating knot groups with the dehn presentation}

\author{Fabienne chouraqui}

\begin{abstract}
Every tame, prime and alternating knot is equivalent to a tame, prime
and alternating knot in regular position, with a common projection.
In this work, we show that the Dehn presentation of the knot group
of a tame, prime, alternating knot, with a regular and common
projection has a finite and complete rewriting system. Although there
are rules in the rewriting system with left-hand side a generator and
which increase the length of the
words we show that the system is terminating.
\end{abstract}
\maketitle

\section{Introduction}

In \cite{2}, we showed that the augmented Dehn presentation of the
knot group of a tame, prime, alternating knot in regular position,
with a common (elementary) projection \cite[p.267]{7} has a finite
and complete rewriting system (with no need of completion) and
this result holds also for alternating links satisfying our
assumptions. We showed there that there are exactly two such
rewriting systems and we gave an algorithm which finds them. This
was carried out using graph theory, applied to the projection of
the knot and to the projection's dual graph. In this paper, we use
one of the complete rewriting systems for the augmented Dehn
presentation of the knot group to find a finite and complete
rewriting system for the Dehn presentation. The same process can
be applied on
the second rewriting system. \\
The main idea is that we can divide the generators
in the Dehn presentation
of the knot group and their inverses  in two disjoint sets, called the sources and the sinks.
This division of the generators is derived from the algorithm which finds
the complete rewriting system
for the augmented Dehn presentation of the knot group and it permits us to compare the words in the free group generated by
these generators in an efficient way.\\ We will describe in some words the way
 we do this and we refer the reader to section $2$ for more details.
If $K$ is a knot and $\wp(K)$ its projection on $R^{2}$,then $\wp(K)$
is said to be \textbf{common} if the boundaries of any two distinct
domains have at most one edge in common and each crossing is on the
boundary of exactly four distinct domains \cite{3}\cite{9}\cite[p267]{7}.
In the \textbf{Dehn presentation} of the knot group (the fundamental
group of $R^{3}\setminus K)$ the generators are
the domains of $R^{2}\setminus \wp(K)$, labelled by $x_{0},x_{1},...,x_{n}$,
starting with the unbounded domain (see Fig.1). The relations arise from
the crossings:
at each crossing, four distinct regions $x_{a}$,$x_{b}$,$x_{c}$,$x_{d}$
meet and the relation arising from this crossing is r:
 $x_{a}\overline{x_{b}}x_{c}\overline{x_{d}}=1$,
 where
$\mathbf{\overline{x}}$ \textbf{denotes the inverse
of $x$}.\\
 The \textbf{Dehn presentation} of the fundamental group
of $R^{3}\setminus$ K is \,\  $<x_{0},x_{1},...,x_{n}\mid r_{1},...,r_{n-1},x_{0}>$
where each $r_{m}$ has the form $x_{i}\overline{x_{j}}x_{k}\overline{x_{l}}=1$
with $x_{i},x_{j},x_{k},x_{l}$ generators and the \textbf{augmented Dehn presentation}
 is the Dehn presentation
with the relation \{$x_{0}=1$\} deleted,
 i.e \,\ $<x_{0},x_{1},...,x_{n}\mid r_{1},...,r_{n-1}>$.

Weinbaum in \cite{9} showed that the augmented Dehn presentation
given by K presents the free product of the knot group
$\pi(R^{3}\setminus K)$
and an infinite cyclic group.\\
The way to divide the generators in the two disjoint sets of sources and sinks
can be described easily as follows:
choose $x_{0}$ to be a source, then all its neighbours will be sinks and so
on iteratively. Next, if $x_{0}$ is a source, then $\overline{x_{0}}$ is
also chosen to be a source and the same holds for all the generators.
This process is consistent,
 since when the knot projection is regular a chess-boarding colouring
of the domains is possible (a regular knot projection is a planar
graph in which every vertex has even degree $4$). The generators
are then renamed according to their being a source or a sink: if
$x_{*}$ is a source then it is renamed $s_{*}$ and if $x_{*}$ is a
sink then it is renamed $t_{*}$ and the relations and the rules
are "rewritten" accordingly. For the eight knot group, if $x_{0}$
and $\overline{x_{0}}$ are sources, then $x_{1}$, $x_{2}$,
$x_{3}$,  $\overline{x_{1}}$, $\overline{x_{2}}$,
$\overline{x_{3}}$ are sinks and $x_{4}$, $x_{5}$,
$\overline{x_{4}}$, $\overline{x_{5}}$ are also sources (see
Fig.1). So,the main result of this paper can be stated as follows:
\begin{thm*}
Let K be a tame, prime and alternating knot with a regular and common
projection. Then the Dehn presentation of the knot group of K has a complete and
finite rewriting system.\\
We shall denote this rewriting system of K by   $\Re''$.
\end{thm*}
\textbf{Example: The augmented Dehn presentation of the eight knot
group }(see Fig.1)

$<s_{0},t_{1},t_{2},t_{3},s_{4},s_{5}\mid
t_{1}\overline{s_{0}}t_{2}\overline{s_{4}},
t_{2}\overline{s_{0}}t_{3}\overline{s_{5}},
t_{1}\overline{s_{5}}t_{3}\overline{s_{0}},
t_{2}\overline{s_{5}}t_{1}\overline{s_{4}}>$

In Section 2, we give some definitions concerning
 rewriting systems. Then,we give
a brief survey of the results obtained in \cite{2} and refer the
reader for more details. In Section 3, we find a rewriting system
$\Re'$ for the Dehn presentation of the knot group and show some
of its properties and then we find a rewriting system $\Re''$
which is equivalent to $\Re'$. In Section 4, we show that the
rewriting system $\Re''$ is complete. The main difficulty is to
show that it is terminating, since there are rules in $\Re''$
which increase the length of the words. \textbf{In this paper, we
shall assume that knots are : tame, prime and alternating with a
regular and common projection. }
\begin{acknowledgment}
I am very grateful to my supervisor Professor Arye Juhasz, for his patience
his encouragement and his very useful and judicious remarks on this paper.
\end{acknowledgment}

\section{Preliminaries}

\subsection{On rewriting systems}
The terminology in this section is from  \cite{1} and  \cite{4}.
Let $\sum$ be a non-empty set. We denote by $\sum^{*}$the free
monoid generated by $\sum$; elements of $\sum^{*}$ are finite
sequences called \textbf{words} and the empty word will be denoted
by 1. The length of the sequence is called the \textbf{length of
the word} $w$ and is denoted by \textbf{$\ell(w)$}.
\begin{defn*}
A \textbf{rewriting system} $\Re$ on $\sum$ is a set of ordered
pairs in $\sum^{*}\times \sum^{*}$ .

If $(l,r) \in \Re$ then for any words $u$ and $v$ in $\sum^{*}$,
 we say the word $ulv$ \textbf{reduces} to the word $urv$ and we
write $ulv\rightarrow urv$ .
\end{defn*}

 A word $w$ is said to be \textbf{reducible}
if there is a word $z$ such that $w\rightarrow z$. If there is no
such $z$ we call $w$ \textbf{irreducible}.

A rewriting system $\Re$ is called \textbf{terminating} \textbf{(or
Noetherian)} if there is no infinite sequence of reductions $w_{1}\rightarrow w_{2}\rightarrow...\rightarrow w_{n}\rightarrow....$

$\Re$ is called \textbf{locally confluent} if for any words $u,v,w$
in $\sum^{*}$ , $w\rightarrow u$ and $w\rightarrow v$ implies that
there is a word $z$ in $\sum^{*}$ such that $u\rightarrow^{*}z$
and $v\rightarrow^{*}z$, where $"\rightarrow^{*}"$ denotes the reflexive
transitive closure of the relation $"\rightarrow"$ .

We call the triple of non-empty words $u,v,w$ in $\sum^{*}$ an \textbf{overlap
ambiguity} if there are $r_{1},r_{2}$ in $\sum^{*}$ such that $uv\rightarrow r_{1}$
and $vw\rightarrow r_{2}$ are rules in $\Re$. We then say that $r_{1}w$
and $ur_{2}$ are the corresponding \textbf{critical pair.} When the triple $u,v,w$ in $\sum^{*}$ is an overlap ambiguity
we will say that the rules $uv\rightarrow r_{1}$ and $vw\rightarrow r_{2}$
\textbf{overlap at $v$} or that there is an \textbf{overlapping}
between the rules $uv\rightarrow r_{1}$ and $vw\rightarrow r_{2}$.
If there exists a word $z$ such that $r_{1}w\rightarrow^{*}z$ and
$ur_{2}\rightarrow^{*}z$ , then we say that the critical pair resulting
from the overlapping of the rules \textbf{}$uv\rightarrow r_{1}$
and $vw\rightarrow r_{2}$ in $\Re$ \textbf{resolves}.

The triple $u,v,w$ of possibly empty words in $\sum^{*}$ is called an \textbf{inclusion
ambiguity} if there are $r_{1},r_{2}$ in $\sum^{*}$(which must be
distinct if both $u$ and $w$ are empty, but otherwise may be equal)
such that $v\rightarrow r_{1}$ and $uvw\rightarrow r_{2}$ are rules
in $\Re$. We then say that $ur_{1}w$ and $r_{2}$ are the corresponding
\textbf{critical pair}. If there exists a word $z$ such
that $ur_{1}w\rightarrow^{*}z$ and $r_{2}\rightarrow^{*}z$ , then
we say that the critical pair resulting from the inclusion ambiguity
of the rule \textbf{}$v\rightarrow r_{1}$ in $uvw\rightarrow r_{2}$
in $\Re$ \textbf{resolves}.

$\Re$ is called \textbf{complete (or convergent)} if $\Re$ is terminating
and locally confluent or in other words if $\Re$ is terminating and
all the critical pairs resolve. So if $\Re$ is complete then every
word $w$ in $\sum^{*}$ has a unique irreducible equivalent word
$z$, which is called the \textbf{normal form} of $w$.

We say that two rewriting systems $\Re$ and $\Re'$ are \textbf{equivalent}
if: $w_{1}\leftrightarrow^{*}w_{2}$ modulo $\Re$ if and only if
$w_{1}\leftrightarrow^{*}w_{2}$ modulo $\Re'$, where $\leftrightarrow^{*}$
denotes the equivalence relation generated by $\rightarrow$ .

We call a rewriting system $\Re$ \textbf{reduced} if for any rule
$l\rightarrow r$ in $\Re$ , $r$ is irreducible and there is no
rule $l'\rightarrow r'$ in $\Re$ such that $l'$ is a subword of
$l$. If $\Re$ is complete then there exists a reduced and complete
rewriting system $\Re$'  which is equivalent to $\Re$ \cite{8}.
The reduced rewriting system $\Re$' which is equivalent to $\Re$
is obtained in two steps: (1) whenever $l\rightarrow r$ is a rule
in $\Re$ and
 $r$ is not irreducible we  reduce  $r$ to its normal form. (2) whenever  $l\rightarrow r$, $l'\rightarrow r'$ are rules in $\Re$ and $l'$ is a subword of
$l$ we  remove  from $\Re$ the rule $l\rightarrow r$.

\subsection{Preliminary results}

Let G be a group presented by $<X\mid R>$, where X is a set of free generators
$\{ x_{0},x_{1},...,x_{n}\}$ and R the set of relations. In order
to define a rewriting system for G, we have to consider the \textbf{monoid
presentation} of G, $<X\cup \overline{X} \mid R\cup R_{0}>$ , where $\overline{X}$ denote the set $\{\overline{x_{0}},...,\overline{x_{n}}\}$
and $R_{0}=\{ x_{i}\overline{x_{i}}=1,\overline{x_{i}}x_{i}=1,$ with  $i\in\{0,1,...,n\}\}$.
The  \textbf{symmetrization process on a relator $r$}  is the following
process: we consider all the relations which can be obtained as consequences
of $r,\overline{r}$ and all their cyclic permutations in a group
and we define $S(r)$ to be the minimal set of relations from which
all the other relations can be derived in the  monoid. Example: if $r:ab=1$
where $a,b\in X\cup\overline{X}$ then we obtain from the symmetrization
process the set of relations $\{ ab=1,ba=1,\overline{b}\,\overline{a}=1,\overline{a}\,\overline{b}=1,b=\overline{a},a=\overline{b}\}$
and $S(r)=\{ a=\overline{b},b=\overline{a}\}$.

We refer the reader to \cite{2} where this process is done in
detail for the relations in the knot group and also to \cite{6}.
In what follows, we will state results from \cite{2} without
proofs.

\begin{rem}
\label{rem:def_S(Rm)}For each relator $r_{m}:x_{i}\overline{x_{j}}x_{k}\overline{x_{l}}=1$,
the set \textbf{S($r_{m}$)} is the following set of relations:
\begin{center}
$\begin{array}{c}
x_{i}\overline{x_{j}}=x_{l}\overline{x_{k}}  \\
\overline{x_{j}}x_{k}=\overline{x_{i}}x_{l}\\
x_{k}\overline{x_{l}}=x_{j}\overline{x_{i}}\\
\overline{x_{l}}x_{i}=\overline{x_{k}}x_{j}\\
\end{array}$
\end{center}
\end{rem}

The set of relations R', which is the union of the sets \{S($r_{m}$)\}
for $1\leq m$ $\leq n-1$ is equivalent to the set of relations $R=\{ r_{m}\mid 1\leq m\leq n-1\}$
for a knot group G, satisfying our assumptions and the problem
is then how to orientate the relations in R' in order to define a
rewriting system.

\textbf{Example: the set of relations R' for the figure-eight knot group}

$\begin{array}{ccccccc}
s_{0}\overline{t_{1}}=t_{2}\overline{s_{4}} &&
s_{0}\overline{t_{2}}=t_{3}\overline{s_{5}} &&
s_{0}\overline{t_{3}}=t_{1}\overline{s_{5}} &&
s_{4}\overline{t_{1}}=t_{2}\overline{s_{5}} \\

\overline{s_{0}}t_{2}=\overline{t_{1}}s_{4}&&
\overline{s_{0}}t_{3}=\overline{t_{2}}s_{5} &&
\overline{s_{0}}t_{1}=\overline{t_{3}}s_{5} &&
\overline{s_{4}}t_{2}=\overline{t_{1}}s_{5}\\

s_{4}\overline{t_{2}}=t_{1}\overline{s_{0}} &&
s_{5}\overline{t_{3}}=t_{2}\overline{s_{0}} &&
s_{5}\overline{t_{1}}=t_{3}\overline{s_{0}} &&
s_{5}\overline{t_{2}}=t_{1}\overline{s_{4}} \\

\overline{s_{4}}t_{1}=\overline{t_{2}}s_{0} &&
\overline{s_{5}}t_{2}=\overline{t_{3}}s_{0} &&
\overline{s_{5}}t_{3}=\overline{t_{1}}s_{0} &&
\overline{s_{5}}t_{1}=\overline{t_{2}}s_{4} \\
\end{array}$\\
Using the fact that all sides of relations in R' have length 2, we
defined a graph $\Delta$, called the derived graph, in the following
way:

\begin{itemize}
\item The vertex-set is the set $X\cup\overline{X}$.
\item For $a,b$ in $X\cup\overline{X}$, there is an oriented edge $a\rightarrow b$
én $\Delta$, if there is a relation in $R'$ such that the word $ab$
is one of its sides.
\end{itemize}
\textbf{Example: The derived  graph  $\Delta$ for the figure-eight knot} (see Fig.2)\\
A subgraph A of $\Delta$ is said to be an \textbf{antipath} in $\Delta$,
if each vertex of A is a source or a sink in A.The number of edges
in this set is called the \textbf{length} of the antipath. Multiple
edges (oriented the same) are not allowed in an antipath. An antipath need not be connected
and an equivalent definition of an antipath is that no oriented path
of length 2 occurs.
An antipath A in $\Delta$ of length half the number of edges in $\Delta$ will define a
rewriting system with no ambiguity overlap between the rules in the following way: all the edges in A
will be the left-hand sides of the rules and all the edges in the complement
of A will be the right-hand sides of rules. \\
\textbf{Example: An antipath in the graph  $\Delta$ for the figure-eight knot} (see Fig.3)\\
It holds that $\Delta$, the derived graph, is the disjoint union
of two Eulerian closed paths of even length in which all the closed
paths have even length. And as such, $\Delta$ admits two disjoint
antipaths, each of length half the number of edges in $\Delta$, which
define a finite and complete rewriting system (with no need of completion).
We will describe the algorithm to find one of the two antipaths $\hat{A}$,
the one we will work with in this paper: choose $x_{0}$ and $\overline{x_{0}}$
to be sources in $\hat{A}$, which means that all the edges going
out of $x_{0}$ and $\overline{x_{0}}$ will belong to $\hat{A}$
and this will determine $\hat{A}$ completely (in an iterative way).
In fact, by having a glance at the knot projection , one can see immediately
for each pair $(x_{i}$,$\overline{x_{i}})$ if these are sources
or sinks in $\hat{A}$ : in the example of the figure-eight knot if $x_{0}$
and $\overline{x_{0}}$ are sources in $\hat{A}$ then $x_{1}$, $\overline{x_{1}}$
,$x_{2}$, $\overline{x_{2}}$, $x_{3}$ and $\overline{x_{3}}$ are
sinks in $\hat{A}$ and $x_{4}, \overline{x_{4}}$, $x_{5}, \overline{x_{5}}$  are also sources
in $\hat{A}$ (by a chess-boarding effect).
We rename the generators according to their being a sink or a source in the antipath $\hat{A}$: if the generator $x_{*}$ is a source it is denoted by $s_{*}$ and if the generator $x_{.}$ is a sink it is denoted by $t_{.}$.
 We will denote by $\Re$
the rewriting system defined by the antipath $\hat{A}$, with the rules "rewritten" accordingly.

\textbf{Example: The rewriting system $\Re$ for the figure-eight knot group:}(see Fig.3)

$\begin{array}{ccccccc}
s_{0}\overline{t_{1}}\rightarrow t_{2}\overline{s_{4}} &&
s_{0}\overline{t_{2}}\rightarrow t_{3}\overline{s_{5}} &&
s_{0}\overline{t_{3}}\rightarrow t_{1}\overline{s_{5}} &&
s_{4}\overline{t_{1}}\rightarrow t_{2}\overline{s_{5}} \\

\overline{s_{0}}t_{2}\rightarrow \overline{t_{1}}s_{4}&&
\overline{s_{0}}t_{3}\rightarrow \overline{t_{2}}s_{5} &&
\overline{s_{0}}t_{1}\rightarrow \overìine{t_{3}}s_{5} &&
\overline{s_{4}}t_{2}\rightarrow \overline{t_{1}}s_{5}\\

s_{4}\overline{t_{2}}\rightarrow t_{1}\overline{s_{0}} &&
s_{5}\overline{t_{3}}\rightarrow t_{2}\overline{s_{0}} &&
s_{5}\overline{t_{1}}\rightarrow t_{3}\overline{s_{0}} &&
s_{5}\overline{t_{2}}\rightarrow t_{1}\overline{s_{4}} \\

\overline{s_{4}}t_{1}\rightarrow \overline{t_{2}}s_{0} &&
\overline{s_{5}}t_{2}\rightarrow \overline{t_{3}}s_{0} &&
\overline{s_{5}}t_{3}\rightarrow \overline{t_{1}}s_{0} &&
\overline{s_{5}}t_{1}\rightarrow \overline{t_{2}}s_{4} \\

s_{i}\overline{s_{i}}\rightarrow 1 &&
 i= 0,4,5 &&
\overline{s_{i}}s_{i}\rightarrow 1&&
  i= 0,4,5\\

t_{i}\overline{t_{i}}\rightarrow 1 &&
 i= 1,2,3 &&
\overline{t_{i}}t_{i}\rightarrow 1&&
  i= 1,2,3\\
\end{array}$\\

\begin{rem}
\label{rem:Each-set-S(r_{m}),}Each set S($r_{m}$), with $r_{m}:x_{i}\overline{x_{j}}x_{k}\overline{x_{l}}$,
gives the following set of rules in $\Re$ if $x_{i},x_{k},\overline{x_{i}},\overline{x_{k}}$
are sources and $x_{j},x_{l},\overline{x_{j}},\overline{x_{l}}$
are sinks in the antipath $\hat{A}$ which defines $\Re$:
\begin{center}
$\begin{array}{c}
(1) \,\,\     s_{i}\overline{t_{j}}\rightarrow t_{l}\overline{s_{k}} \\
(2)\,\,\      \overline{s_{i}}t_{l}\rightarrow\overline{t_{j}}s_{k}\\
(3) \,\,\     s_{k}\overline{t_{l}}\rightarrow t_{j}\overline{s_{i}}\\
(4)\,\,\      \overline{s_{k}}t_{j}\rightarrow\overline{t_{l}}s_{i}\\
\end{array}$
\end{center}
If $x_{j},x_{l},\overline{x_{j}},\overline{x_{l}}$ are sources  and $x_{i},x_{k},\overline{x_{i}},\overline{x_{k}}$ are sinks in  the antipath $\hat{A}$, then a set of rules of the same kind is obtained. Note that since $x_{0}$ and $\overline{x_{0}}$
are sources in $\hat{A}$, $x_{0}$ and $\overline{x_{0}}$ will appear
as the first letter of the left-hand side of a rule or as the last
letter of the right-hand side.
\end{rem}

\section{definition of a rewriting system $\Re'$ for the dehn presentation }

\subsection{Definition of $\Re'$ }

We will use in what follows some tools developped  in  \cite{9}.
Let K be a tame, prime and alternating knot whose projection is
regular and common. Let us denote the knot group
$\pi(R^{3}\setminus K)$ by G and assume its Dehn presentation is
$<x_{0},x_{1},...,x_{n}\mid r_{1},...,r_{n-1},x_{0}>$. So H , the
group presented by the augmented Dehn presentation of $\pi(R^{3}
\setminus  K)$ is $<x_{0},x_{1},...,x_{n}\mid r_{1},...,r_{n-1}>$.
Let $\Re$ be the complete and finite rewriting system for H
obtained in the previous section. Let F be the free group on n+1
generators  $x_{0},x_{1},...,x_{n}$. Let $\Phi$ be the
endomorphism of F determined by: $\Phi(x_{0})=1$
and $\Phi(x_{j})=x_{j}, 1\leq j\leq n$. \\
Then $G\simeq<x_{1},...,x_{n}\mid\Phi(r_{1}),...,\Phi(r_{n-1})>$.

In order to define a rewriting system $\Re'$ for G, we need to apply
the symmetrization process on $\Phi(r_{m})$ , for $1\leq m\leq n-1$
and to add of course the relations $x_{i}\overline{x_{i}}=1,\overline{x_{i}}x_{i}=1$
for $1\leq i\leq n$.

We recall that each relator has the form:
$r_{m}:x_{i}\overline{x_{j}}x_{k}\overline{x_{l}}$, so if none of
the indices $i,j,k,l$ is 0 then
$\Phi(r_{m})=x_{i}\overline{x_{j}}x_{k}\overline{x_{l}}$ and the
symmetrization process applied on $\Phi(r_{m})$ gives the set of
relations S($r_{m}$) (see remark \ref{rem:def_S(Rm)}). If one of
the indices is 0, assume $i=0$, then
$\Phi(r_{m})=\overline{x_{j}}x_{k}\overline{x_{l}}$ and the
symmetrization process applied on $\Phi(r_{m})$ gives the
following set of relations, denoted by $S(\Phi(r_{m}))$:
$\begin{array}{ccc} (a) x_{l}=\overline{x_{j}}x_{k} &&
(d)\overline{x_{l}}=\overline{x_{k}}x_{j}\\
(b) x_{j}=x_{k}\overline{x_{l}}
&&
(e) \overline{x_{j}}=x_{l}\overline{x_{k}}\\
(c)\overline{x_{k}}=\overline{x_{l}}\,\overline{x_{j}}
&&
(f) x_{k}=x_{j}x_{l}\\
\end{array}$

\textbf{We denote the relations
$x_{i}\overline{x_{i}}=1,\overline{x_{i}}x_{i}=1$ for $1\leq i\leq
n$ by (0), and these relations will belong to all the sets of
relations we will consider, even if we don't mention this
explicitly}.
We will work all along with the assumption that $i=0$, where $i$ is the index of the first letter in $\Phi(r_{m})$. Clearly, there is no loss of generality in doing this.\\

Now, let reverse the order of the operations on a relation $r_{m}$,
i.e. first making the symmetrization process to obtain $S(r_{m})$and
then applying $\Phi$ on both sides of the relations in $S(r_{m})$.The
set of relations obtained will be denoted by $\Phi(S(r_{m}))$. If
none of the indices $i,j,k,l$ in $r_{m}$ is 0 then $\Phi(S(r_{m}))=S(\Phi(r_{m}))=S(r_{m})$
and if $i=0$ then $\Phi(S((r_{m})))$ is the following set of relations:
$\begin{array}{cc}
x_{l}=\overline{x_{j}}x_{k}&
\overline{x_{l}}=\overline{x_{k}}x_{j} \\
x_{j}=x_{k}\overline{x_{l}} &
\overline{x_{j}}=x_{l}\overline{x_{k}} \\
\end{array}$
 which correspond respectively to the relations $(a),(b),(d),(e)$
in $S(\Phi(r_{m}))$.\\
 So, by reversing the order of the operations
we lost the two relations
 $(c): \overline{x_{k}}=\overline{x_{l}}\,\overline{x_{j}}$ and
$(f): x_{k}=x_{j}x_{l}$.
Yet, we will show that the equivalence relation generated by $\Phi(S(r_{m}))$
is the same as that generated by $S(\Phi(r_{m}))$ and that in fact
it is sufficient to consider an even smaller set of relations.

\begin{claim}
The equivalence relation generated by the following set of relations,
denoted by $\Phi_{m} :
\begin{array}{c}
(1)x_{l}=\overline{x_{j}}x_{k}\\
(2)x_{j}=x_{k}\overline{x_{l}} \\
(3)\overline{x_{k}}=\overline{x_{l}}\,\overline{x_{j}} \\
\end{array}$ is the same as that generated by $\Phi(S(r_{m}))$ and $S(\Phi(r_{m}))$.
\end{claim}
\begin{proof}
First, the relations $(d),(e)$ and $(f)$ can be derived from the relations
$(1),(2),(3),(0)$ in the following way: \\
$(d):\overline{x_{k}}x_{j}=\overline{x_{k}}(x_{k}\overline{x_{l}})=\overline{x_{l}}$
  by using $(1)$and then $(0)$\\
$(e):x_{l}\overline{x_{k}}=(\overline{x_{j}}x_{k})\overline{x_{k}}=\overline{x_{j}}$
  by using $(2)$ and then $(0)$\\
 $(f):x_{j}x_{l}=(x_{k}\overline{x_{l}})(\overline{x_{j}}x_{k})=x_{k}(\overline{x_{l}}\,\overline{x_{j}})x_{k}=x_{k}(\overline{x_{k}})x_{k}=x_{k}$
  by using $(1),(2),(3),(0)$.\\
 So, the equivalence relation generated
by $S(\Phi(r_{m}))$ is the same as that generated by $\Phi_{m}$.\\

It remains to show that the relation $(3)$can be derived from the
relations in $\Phi(S(r_{m}))$: $(3)$ is obtained from $(d)$ or
$(e)$ in one of the following ways: from $(d)$:
 $\overline{x_{l}}\,\overline{x_{j}}=(\overline{x_{k}}x_{j})\overline{x_{j}}=\overline{x_{k}}$
or from $(e)$:
$\overline{x_{l}}\,\overline{x_{j}}=\overline{x_{l}}(x_{l}\overline{x_{k}})=\overline{x_{k}}$.
So, the equivalence relation generated by $\Phi(S(r_{m}))$ is the
same as that generated by $\Phi_{m}$.
\end{proof}
In order to define a rewriting system $\Re'$ for G, we have to
orientate the relations in $\Phi_{m}$, for $1\leq m\leq n-1$.
Here, the complete rewriting system $\Re$ found for H plays an
important role: we will orientate the relations in $\Phi_{m}$ in
the same way as the corresponding relations in $S(r_{m})$. When
$x_{0}$ does not appear in the relation $r_{m}$ then things are
clear since $\Phi_{m}=S(r_{m})$ but when $x_{0}$ appears in
$r_{m}$ then we have to recover for each relation in $\Phi_{m}$
what is the corresponding relation in $S(r_{m})$. It seems to be a
hard task but in fact it is not, since we know exactly how rules
look like in $\Re$ and these are described in remark
\ref{rem:Each-set-S(r_{m}),}. So, one can see immediately that the
relation $(3)$ in $\Phi_{m}$ does not correspond to any relation
in $S(r_{m})$ and one can guess that the relations corresponding
to $(1)$ and $(2)$ in $S(r_{m})$ are respectively
$\overline{x_{0}}x_{l}=\overline{x_{j}}x_{k}$ and
$x_{j}\overline{x_{0}}=x_{k}\overline{x_{l}}$. Since
$x_{0},\overline{x_{0}}$ are sources in the antipath $\hat{A}$
which defines $\Re$, we have that $x_{k},\overline{x_{k}}$ are
also sources and $x_{j},\overline{x_{j}},x_{l},\overline{x_{l}}$
are sinks in  $\hat{A}$. So, in $\Re$ these relations are
orientated $\overline{s_{0}}t_{l} \rightarrow
\overline{t_{j}}s_{k}$ and $s_{k}\overline{t_{l}} \rightarrow
t_{j}\overline{s_{0}}$ respectively (see remark
\ref{rem:Each-set-S(r_{m}),}). Since we want to orientate the
relations in $\Phi_{m}$ in the same way as the corresponding
relations in $S(r_{m})$, we orientate the relations $(1)$ and
$(2)$
 in the following way:\\
$(1):t_{l} \rightarrow\overline {t_{j}}s_{k}$ \\
$(2):s_{k}\overline{t_{l}}\rightarrow t_{j}$\\
The following orientation is chosen for the relation $(3)$ in
$\Phi_{m}$:\\
$(3): \overline{t_{l}}\,\overline{t_{j}}\rightarrow\overline{s_{k}}$\\
 We will denote this set of 3 rules by $\vec{\Phi_{m}}$. \\

The rewriting system $\Re'$ for $G$ is the union of the following sets of rules:

\begin{tabular}{|c||c||c|}
\hline
 none of the indices 0&
 one of the indices 0 &
\\
\hline

$ s_{i}\overline{t_{j}}\rightarrow t_{l}\overline{s_{k}} $&
$(1):t_{l}\rightarrow\overline{t_{j}}s_{k}$ &
$(0):s_{i}\overline{s_{i}}\rightarrow1$\\
$ \overline{s_{i}}t_{l}\rightarrow\overline{t_{j}}s_{k}$&
$(2):s_{k}\overline{t_{l}}\rightarrow t_{j}$ &
$(0):\overline{s_{i}}s_{i}\rightarrow1$\\
$ s_{k}\overline{t_{l}}\rightarrow t_{j}\overline{s_{i}}$&
$(3):\overline{t_{l}}\,\overline{t_{j}}\rightarrow\overline{s_{k}}$ &
$(0):t_{i}\overline{t_{i}}\rightarrow1$\\
$  \overline{s_{k}}t_{j}\rightarrow\overline{t_{l}}s_{i}$&&
$(0):\overline{t_{i}}t_{i}\rightarrow1$\\
\hline
\end{tabular}

\begin{rem}\label{rem:t+_neighX0}
There is a rule in $\Re'$ whose left-hand side is
a  sink $t$ with positive exponent if and only if this $t$ is
connected by an edge to $\overline{x_{0}}$ in $\hat{A}$, i.e it corresponds
 to a neighbouring region of the unbounded domain in the knot projection.
\end{rem}

\textbf{Example: The rewriting system $\Re'$ for the eight knot group}
$\begin{array}{ccccccc}
t_{2}\rightarrow\overline{t_{1}}s_{4} &&
t_{3}\rightarrow\overline{t_{2}}s_{5} &&
t_{1}\rightarrow\overline{t_{3}}s_{5} &&
s_{4}\overline{t_{1}}\rightarrow
t_{2}\overline{s_{5}} \\

s_{4}\overline{t_{2}}\rightarrow t_{1} &&
s_{5}\overline{t_{3}}\rightarrow t_{2} &&
s_{5}\overline{t_{1}}\rightarrow t_{3} &&
\overline{s_{4}}t_{2}\rightarrow\overline{t_{1}}s_{5}\\

\overline{t_{2}}\,\overline{t_{1}}\rightarrow\overline{s_{4}} &&
\overline{t_{3}}\,\overline{t_{2}}\rightarrow\overline{s_{5}} &&
\overline{t_{1}}\,\overline{t_{3}}\rightarrow\overline{s_{5}} &&
s_{5}\overline{t_{2}}\rightarrow t_{1}\overline{s_{4}}\\

&&
 &&
 &&
\overline{s_{5}}t_{1}\rightarrow\overline{t_{2}}s_{4}\\

s_{i}\overline{s_{i}}\rightarrow 1 &&
 i= 0,4,5 &&
\overline{s_{i}}s_{i}\rightarrow 1&&
  i= 0,4,5\\

t_{i}\overline{t_{i}}\rightarrow 1 &&
 i= 1,2,3 &&
\overline{t_{i}}t_{i}\rightarrow 1&&
  i= 1,2,3\\
\end{array}$ \\
 Note that $\Re'$ is not complete. As an example, the word $\overline{t_{3}}\,\overline{t_{2}}\,\overline{t_{1}}$ reduces to two different irreducible words $\overline{s_{5}}\,\overline{t_{1}}$ and $\overline{t_{3}}\,\overline{s_{4}}$.

\subsection{Properties of $\Re'$ }

\begin{claim}
\label{cla:rel_in_altern}If there is a relator $r_{m}:x_{i}\overline{x_{j}}x_{k}\overline{x_{l}}$
in the augmented Dehn presentation of an alternating knot group, then
there are $1 \leq k',l' \leq n$ such that there is  also the relator $r_{m'}:x_{k'}\overline{x_{j}}x_{i}\overline{x_{l'}}$.
\end{claim}
\begin{proof}
Since the knot is alternating and we have $r_{m}:x_{i}\overline{x_{j}} x_{k} \overline{x_{l}}$,
then we have necessarily in the knot projection the  situation described  in Fig.4.
This  means that there is also the relation $r_{m'}:x_{k'} \overline{x_{j}} x_{i} \overline{x_{l'}}$
(or its cyclic permutation $r_{m'}:x_{i} \overline{x_{l'}} x_{k'} \overline{x_{j}}$)
\end{proof}
\begin{claim}
\label{cla:if_setofm_then_setof m'}If the set $\vec{\Phi_{m}}$:
$\begin{array}{c}
(1_{l,j,k}):t_{l}\rightarrow\overline{t_{j}}s_{k} \\
(2_{l,j,k}):s_{k}\overline{t_{l}}\rightarrow t_{j} \\
(3_{l,j,k}):\overline{t_{l}}\,\overline{t_{j}}\rightarrow\overline{s_{k}} \\
\end{array}$ is in $\Re'$ then\\
 $\vec{\Phi_{m'}}: \begin{array}{c}
(1_{j,l',k'}):t_{j}\rightarrow\overline{t_{l'}}s_{k'}\\
(2_{j,l',k'}):s_{k'}\overline{t_{j}}\rightarrow t_{l'} \\
(3_{j,l',k'}):\overline{t_{j}}\,\overline{t_{l'}}\rightarrow\overline{s_{k'}} \\
\end{array}$ is also in $\Re'$.
\end{claim}
\begin{proof}
It follows from claim \ref{cla:rel_in_altern} , where $\vec{\Phi_{m}}$
and $\vec{\Phi_{m'}}$ are the sets of rules in $\Re'$ corresponding to
$r_{m}$ and $r_{m'}$ respectively.
\end{proof}

\subsection{Definition of an ordering on the words in $(X\cup\overline{X})^{*}$}
By abuse of notation, we denote also by $X\cup\overline{X}$ the union of the set
 of sinks and the set of sources, i.e the set of generators
and their inverses after their renaming.
We write a non-empty word $w$ in $(X\cup\overline{X})^{*}$,
$w=S_{1}t_{1}S_{2}t_{2}..S_{k}t_{k}S_{k+1}$,
with $S_{i}$ a sequence of sources $s_{1}s_{2}..s_{n_{i}}$
such that the length of $S_{i}$ (i.e the number of occurrences of $s$), denoted by $n_{i}$,
satisfies $0\leq n_{i} \leq \ell(w)$ for $1\leq i\leq k+1$.

We will define an ordering on the words in $(X\cup\overline{X})^{*}$
in the following way.
Let $>$ denote the usual ordering on the nonnegative integers and
let $>_{lex}$ denote the lexicographic ordering from the left on $4-$tuples
of nonnegative integers, induced by $>$. Define \\

$V_{1}(w)$ to be the number of occurrences of sinks with positive and negative
 exponent (the number of $t$) in $w$. It holds that $V_{1}(w) \leq \ell(w)$.\\

$V_{2}(w)$ to be the number of occurrences of sinks with positive exponent in $w$, but only
 those sinks which are connected by an edge to $\overline{x_{0}}$ in the antipath $\hat{A}$ which defines $\Re$. These generators correspond in fact to neighbouring regions of the unbounded region ($x_{0})$
in the knot projection. It holds that $V_{2}(w) \leq \ell(w)$.\\

$V_{3}(w)=\sum _{j=1}^{j=k}\sum_{i=0}^{i=k-j}2^{i}n_{j}$ \\
It holds that $V_{3}(w) \leq 2^{2k} \ell(w)$.

$V_{4}(w)$ to be $\ell(w)$, the length of  $w$.\\

$V(w)=(V_{1}(w), V_{2}(w), V_{3}(w), V_{4}(w))$ \\
Observe that $V(w) \geq (0,0,0,0)$
for every non-empty word $w$.
Define $w > w'$ if and only if $V(w) >_{lex} V(w')$ \\

\textbf{Example: Computation of $V(t_{2})$ and
$V(\overline{t_{1}}s_{4})$} It holds that $V_{1}(t_{2}) =
V_{1}(\overline{t_{1}}s_{4}) = 1$ but $V_{2}(t_{2}) = 1 $ and
$V_{2}(\overline{t_{1}}s_{4}) = 0$. So, $V(t_{2}) >_{lex}
V(\overline{t_{1}}s_{4}) $.

\section{Definition of an equivalent rewriting system $\Re''$}

\begin{lem}
\label{lem:add_rule_super_3}Assume $\vec{\Phi_{m}}:
\begin{array}{c}
(1_{l,j,k}):t_{l}\rightarrow\overline{t_{j}}s_{k}\\
(2_{l,j,k}):s_{k}\overline{t_{l}}\rightarrow t_{j} \\
(3_{l,j,k}):\overline{t_{l}}\,\overline{t_{j}}\rightarrow\overline{s_{k}} \\
\end{array}$ and $\vec{\Phi_{m'}} :
\begin{array}{c}
(1_{j,l',k'}):t_{j}\rightarrow\overline{t_{l'}}s_{k'}\\
(2_{j,l',k'}):s_{k'}\overline{t_{j}}\rightarrow t_{l'} \\
(3_{j,l',k'}):\overline{t_{j}}\,\overline{t_{l'}}\rightarrow\overline{s_{k'}} \\
\end{array}$ are in $\Re'$. Then the rule
$(4)\overline{s_{k}}\,\overline{t_{l'}}\rightarrow\overline{t_{l}}\,\overline{s_{k'}}$
is obtained from the equivalence relation
generated by rules of kind $(3)$.
\end{lem}
\begin{proof}
We have
$\overline{s_{k}}\,\overline{t_{l'}}=\overline{t_{l}}\,\overline{t_{j}}\,\overline{t_{l'}}=
\overline{t_{l}}\,\overline{s_{k'}}$,
using first the equivalence relation generated by $(3_{l,j,k})$
and then by $(3_{j,l',k'})$. In order to orientate this relation,
we use the ordering defined above: $V(\overline{s_{k}}\,\overline{t_{l'}})=(1,0,1,2)$ and
$V(\overline{t_{l}}\,\overline{s_{k'}})=(1,0,0,2)$.
So we have $\overline{s_{k}}\,\overline{t_{l'}}\rightarrow\overline{t_{l}}\,\overline{s_{k'}}$.
\end{proof}
\begin{rem}
Adding the rule
$(4)\overline{s_{k}}\,\overline{t_{l'}}\rightarrow\overline{t_{l}}\,\overline{s_{k'}}$
to $\Re'$ will resolve the critical pairs resulting from the
overlapping of the rules $(3_{l,j,k})$ and $(3_{j,l',k'})$.
\end{rem}

\textbf{Example:} In the rewriting system $\Re'$ for the figure-eight knot group we have the following two rules of kind (3):
$\overline{x_{2}}\,\overline{x_{1}}\rightarrow\overline{x_{4}}$  and
$\overline{x_{1}}\,\overline{x_{3}}\rightarrow\overline{x_{5}}$.
The rule  of kind $(4)$ obtained  from the equivalence relation generated by these two rules is then:
$\overline{x_{4}}\,\overline{x_{3}}\rightarrow \overline{x_{2}}\,\overline{x_{5}}$.

\begin{lem}
\label{lem:Add_rule_incl_over}Assume $\vec{\Phi_{m}}: \begin{array}{c}
(1_{l,j,k}):t_{l}\rightarrow\overline{t_{j}}s_{k} \\
(2_{l,j,k}):s_{k}\overline{t_{l}}\rightarrow t_{j} \\
(3_{l,j,k}):\overline{t_{l}}\,\overline{t_{j}}\rightarrow\overline{s_{k}} \\
\end{array}$ and
$\begin{array}{c}
(\alpha) s_{m}\overline{t_{p}}\rightarrow t_{l}\overline{s_{m'}} \\
(\beta) s_{m'}\overline{t_{l}} \rightarrow t_{p}\overline{s_{m}}\\
(\gamma) \overline{s_{m}}t_{l} \rightarrow  \overline{t_{p}}s_{m'}\\
(\delta) \overline{s_{m'}}t_{p}\rightarrow\overline{t_{l}}s_{m}\\
\end{array}$ are in $\Re'$.
Then the rule
 $(5)\overline{s_{m}}\,\overline{t_{j}}\rightarrow\overline{t_{p}}s_{m'}\,\overline{s_{k}}$
is obtained from
the equivalence relation generated by $(3_{l,j,k})$ and $(\beta)$.
Further, rule $(\gamma)$ can be derived using the rule $(5)$.
\end{lem}
\begin{proof}
We have
$\overline{t_{p}}s_{m'}\,\overline{s_{k}}=\overline{t_{p}}s_{m'}\overline{t_{l}}\,
\overline{t_{j}}=\overline{t_{p}}t_{p}\overline{s_{m}}\,\overline{t_{j}}=\overline{s_{m}}\,\overline{t_{j}}$
, using first the equivalence relation generated by $(3_{l,j,k})$
and then by $(\beta)$ and $(0)$. In order to orientate this
relation, we use the ordering defined above:
$V(\overline{s_{m}}\,\overline{t_{j}})=(1,0,1,2)$ and
$V(\overline{t_{p}}s_{m'}\,\overline{s_{k}} )=(1,0,0,3)$. So, we
have
$\overline{s_{m}}\,\overline{t_{j}}\rightarrow\overline{t_{p}}s_{m'}\,\overline{s_{k}}$.
 Further, we have
$\overline{s_{m}}t_{l}\rightarrow
\overline{s_{m}}\,\overline{t_{j}}s_{k}\rightarrow
\overline{t_{p}}s_{m'}\overline{s_{k}}s_{k}
\rightarrow\overline{t_{p}}s_{m'}$
by using rules $(1_{l,j,k})$, $(5)$ and (0).
\end{proof}
\begin{rem}\label{rem:no_t+}
Adding the rule
$(5)\overline{s_{m}}\,\overline{t_{j}}\rightarrow\overline{t_{p}}s_{m'}\,\overline{s_{k}}$
to $\Re'$ will resolve the critical pairs resulting from the overlapping of the rules
$(3_{l,j,k})$ and $\beta$ and from the inclusion ambiguity of rule $(1_{l,j,k})$ and $\gamma$.
Additionally, removing the rule $\gamma$ from $\Re'$ and replacing it by $(5)$ implies that
$t_{l}$ cannot appear in the left-hand side of any rule except in $(1_{l,j,k})$.
\end{rem}

\textbf{Example:} In the rewriting system $\Re'$ for the figure-eight knot group we have
the two following rules of kind (3) and
($\beta$) respectively:
$\overline{x_{2}}\,\overline{x_{1}}\rightarrow\overline{x_{4}}$ and
$x_{5}\overline{x_{2}}\rightarrow x_{1}\overline{x_{4}}$.
The  rule of kind $(5)$ obtained  from the equivalence relation generated by
these two rules is then:
 $\overline{x_{4}}\,\overline{x_{1}}\rightarrow\overline{x_{1}}x_{5}\overline{x_{4}}$.

\begin{prop}\label{prop:R''_reduced}
\label{pro: R''_equival_R'}Assume that  $\Re'$  is the union of the following sets of rules:

\begin{tabular}{|c||c||c|}
\hline
&
$\vec{\phi_{m}}$&
$\vec{\phi_{m'}}$
\\
\hline

$(\alpha) s_{m}\overline{t_{p}}\rightarrow t_{l}\overline{s_{m'}} $&
$(1_{l,j,k}):t_{l}\rightarrow\overline{t_{j}}s_{k}$ &
$(1_{j,l',k'}):t_{j}\rightarrow\overline{t_{l'}}s_{k'}$\\

$(\beta) s_{m'}\overline{t_{l}}\rightarrow t_{p}\overline{s_{m}}$&
$(2_{l,j,k}):s_{k}\overline{t_{l}}\rightarrow t_{j}$ &
$(2_{j,l',k'}):s_{k'}\overline{t_{j}}\rightarrow t_{l'}$\\

$(\gamma)\overline{s_{m}}t_{l}\rightarrow \overline{t_{p}}s_{m'}$&
$(3_{l,j,k}):\overline{t_{l}}\,\overline{t_{j}}\rightarrow\overline{s_{k}}$ &
$(3_{j,l',k'}):\overline{t_{j}}\,\overline{t_{l'}}\rightarrow\overline{s_{k'}}$\\

$(\delta)\overline{s_{m'}}t_{p}\rightarrow\overline{t_{l}}s_{m}$ &&
\\
\hline
\end{tabular}\\
Let $\Re''$ be the rewriting system obtained by adding to $\Re'$
the rules $(4)\overline{s_{k}}\,\overline{t_{l'}}\rightarrow\overline{t_{l}}\,\overline{s_{k'}}$
and $(5)\overline{s_{m}}\,\overline{t_{j}}\rightarrow\overline{t_{p}}s_{m'}\,\overline{s_{k}}$,
replacing $(\alpha)$ by
$(\alpha'):s_{m}\overline{t_{p}}\rightarrow\overline{t_{j}}s_{k}\overline{s_{m'}}$ and $(2_{j,l,k})$ by
$(2'_{j,l,k}):s_{k}\overline{t_{l}}\rightarrow\overline{t_{l'}}s_{k'}$ and removing $(\gamma)$.
That means that  $\Re''$ is the union of the following sets of rules:

\begin{tabular}{|c||c||c|}
\hline
&$\vec{\phi_{m}}$&
$\vec{\phi_{m'}}$
\\
\hline

$(\alpha') s_{m}\overline{t_{p}}\rightarrow \overline{t_{j}}s_{k}\overline{s_{m'}} $&
$(1_{l,j,k}):t_{l}\rightarrow\overline{t_{j}}s_{k}$ &
$(1_{j,l',k'}):t_{j}\rightarrow\overline{t_{l'}}s_{k'}$\\

$(\beta) s_{m'}\overline{t_{l}}\rightarrow t_{p}\overline{s_{m}}$&
$(2'_{l,j,k}):s_{k}\overline{t_{l}}\rightarrow \overline{t_{l'}}s_{k'}$ &
$(2_{j,l',k'}):s_{k'}\overline{t_{j}}\rightarrow t_{l'}$\\

$(\delta)\overline{s_{m'}}t_{p}\rightarrow\overline{t_{l}}s_{m}$ &
$(3_{l,j,k}):\overline{t_{l}}\,\overline{t_{j}}\rightarrow\overline{s_{k}}$ &
$(3_{j,l',k'}):\overline{t_{j}}\,\overline{t_{l'}}\rightarrow\overline{s_{k'}}$\\
\hline
$(5):\overline{s_{m}}\,\overline{t_{j}}\rightarrow\overline{t_{p}}s_{m'}\,\overline{s_{k}}$&
$(4):\overline{s_{k}}\,\overline{t_{l'}}\rightarrow\overline{t_{l}}\,\overline{s_{k'}}$&
 \\
\hline
\end{tabular}

Then the  rewriting system
$\Re''$ is equivalent to $\Re'$ and is reduced.
\end{prop}
\begin{proof}
From lemmas \ref{lem:add_rule_super_3} and \ref{lem:Add_rule_incl_over},
the rules $(4)$ and $(5)$ are obtained from the equivalence relation
generated by the rules in $\Re'$. $(\alpha')$ and $(2'_{j,l,k})$
are obtained by reducing the right-hand side of $(\alpha)$ and $(2_{j,l,k})$
using the rule $(1_{l,j,k})$ and $(\gamma)$ can be derived
using $(5)$. So, $\Re''$ is equivalent to $\Re'$. $\Re''$ is reduced
since the right-hand sides of the rules in $\Re''$ have been reduced
and the rules (of kind $\gamma$) for which the left-hand side contains a subword which
is the left-hand side of an other rule have been removed.
\end{proof}
\begin{rem}
If additionally there is a rule of kind (1) which has $t_{p}$ as
left-hand side then the rule $(\delta)$ is replaced by a rule of
kind $(5)$ and we replace $(\beta)$ by $(\beta')$ in which the
right-hand side is reduced.
\end{rem}
\textbf{Example: The rewriting system $\Re''$ for the figure-eight knot
group }

$\begin{array}{ccccccc}
t_{2}\rightarrow\overline{t_{1}}s_{4} &&
t_{3}\rightarrow\overline{t_{2}}s_{5} &&
t_{1}\rightarrow\overline{t_{3}}s_{5} &&
s_{4}\overline{t_{1}}\rightarrow
\overline{t_{1}}s_{4}\overline{s_{5}} \\

s_{4}\overline{t_{2}}\rightarrow \overline{t_{3}}s_{5}  &&
s_{5}\overline{t_{3}}\rightarrow \overline{t_{1}}s_{4} &&
s_{5}\overline{t_{1}}\rightarrow \overline{t_{2}}s_{5} &&
\overline{s_{4}}\overline{t_{1}}\rightarrow\overline{t_{1}}s_{5}\overline{s_{4}}\\

\overline{t_{2}}\,\overline{t_{1}}\rightarrow\overline{s_{4}} &&
\overline{t_{3}}\,\overline{t_{2}}\rightarrow\overline{s_{5}} &&
\overline{t_{1}}\,\overline{t_{3}}\rightarrow\overline{s_{5}} &&
s_{5}\overline{t_{2}}\rightarrow \overline{t_{3}}s_{5}\overline{s_{4}}\\

\overline{s_{4}}\,\overline{t_{3}}\rightarrow\overline{t_{2}}\,\overline{s_{5}}&&
\overline{s_{5}}\,\overline{t_{2}}\rightarrow\overline{t_{1}}\,\overline{s_{5}}&&
\overline{s_{5}}\,\overline{t_{1}}\rightarrow\overline{t_{3}}\,\overline{s_{4}}&&
\overline{s_{5}}\,\overline{t_{3}}\rightarrow\overline{t_{2}}s_{4}\overline{s_{5}}\\

s_{i}\overline{s_{i}}\rightarrow 1 &&
 i= 0,4,5 &&
\overline{s_{i}}s_{i}\rightarrow 1&&
  i= 0,4,5\\

t_{i}\overline{t_{i}}\rightarrow 1 &&
 i= 1,2,3 &&
\overline{t_{i}}t_{i}\rightarrow 1&&
  i= 1,2,3\\

\end{array}$\\

Note that $\overline{s_{5}}t_{1}\rightarrow\overline{t_{2}}s_{4}$
and $\overline{s_{4}}t_{2}\rightarrow\overline{t_{1}}s_{5}$ (in
the fourth set of rules) have been replaced by the rules of kind $(5)$,
$\overline{s_{5}}\,\overline{t_{3}}\rightarrow\overline{t_{2}}s_{4}\overline{s_{5}}$
and $\overline{s_{4}}\,\overline{t_{1}}\rightarrow\overline{t_{1}}s_{5}\overline{s_{4}}$
respectively.
The rewriting system $\Re''$ for the trefoil knot
group  is: $\Re''=\{x_{1}\rightarrow\overline{x_{4}}x_{3},
x_{3}\overline{x_{1}}\rightarrow\overline{x_{2}}x_{3}, \overline{x_{1}}\,\overline{x_{4}}\rightarrow
\overline{x_{3}}\}\cup\{x_{2}\rightarrow\overline{x_{1}}x_{3},
x_{3}\overline{x_{2}}\rightarrow\overline{x_{4}}x_{3}, \overline{x_{2}}\,\overline{x_{1}}\rightarrow
\overline{x_{3}}\} \cup\{ x_{4}\rightarrow\overline{x_{2}}x_{3},
x_{3}\overline{x_{4}}\rightarrow\overline{x_{1}}x_{3},
\overline{x_{4}}\,\overline{x_{2}}\rightarrow\overline{x_{3}}\}
\cup\{ \overline{x_{3}}\,\overline{x_{2}}\rightarrow\overline{x_{1}}\,\overline{x_{3}},
\overline{x_{3}}\,\overline{x_{1}}\rightarrow\overline{x_{4}}\,\overline{x_{3}},
\overline{x_{3}}\,\overline{x_{4}}\rightarrow\overline{x_{2}}\,\overline{x_{3}}\}\cup
\{x_{i}\overline{x_{i}}\rightarrow1, \overline{x_{i}}x_{i}\rightarrow1, for 0\leq i\leq4\}$
(see \cite{2} for notation) Note that for the trefoil knot  there are no rules of kind $\alpha,\beta,\gamma,\delta$ and
$(5)$ since all the regions have an edge in common with
the unbounded region $x_{0}$.\\

The rewriting system $\Re''$ is complete and this is what we will
show in a more general context in the following.

\section{completeness of the rewriting system $\Re''$}

\subsection{Termination of the rewriting system $\Re''$ }

As before, we will use the following notation: $t$ denotes a sink in the antipath $\hat{A}$
with positive or negative exponent, $t^{+}$ denotes a sink
with positive exponent which is connected by an edge to $\overline{x_{0}}$ in the antipath $\hat{A}$. $s$ denotes a source in $\hat{A}$
with positive or negative exponent and
$S_{i}$ a sequence $s_{1}s_{2}..s_{n_{i}}$.

We will show that we can divide schematically all the rules in $\Re''$ (except
rules of kind $(0)$) in
the following four classes:
$\begin{array}{c}
(A)\,\,\ st\rightarrow tss\\
(B)\,\,\ st\rightarrow ts\\
(C)\,\,\ t^{+}\rightarrow ts\\
(D)\,\,\ tt\rightarrow s\\
\end{array}$ \\
Note that the $s,t$ occurring in one side of a relation are different
from the $s,t$ in the other side, although we use the same letters.
\begin{lem}
Each rule in $\Re''$ belongs to one of the classes (A),(B),(C) or (D) described above.
\end{lem}
\begin{proof}
The rules of kind $(5), \alpha'$ and $\beta'$ belong to the class (A). The rules of kind $(2'),(4)$ and $\alpha , \beta ,  \gamma, \delta$  belong to class (B).
The rules of kind $(1)$ belong to class (C)and the
 rules of kind $(3)$ belong to class (D).
\end{proof}

\begin{lem}\label{rem:t+_onlyC}
In $\Re''$, if a specific $t^{+}$ occurs in a rule from class (C),
then the same $t^{+}$ cannot occur in any other rule, except in rules of kind
$(0)$.
\end{lem}
\begin{proof}
From  proposition \ref{prop:R''_reduced}, we have that $\Re''$ is reduced, so a
$t^{+}$ which occurs in a rule from class (C) cannot occur in the right-hand side of any rule.
From the description of the rules in $\Re'$  and $\Re''$ in proposition \ref{prop:R''_reduced}, we have that a  $t$ with positive exponent sign can occur only in rules
 of kind $(1)$, $\gamma$, $\delta$ and $(0)$, i.e in rules from class (C) or (B) respectively. Now, if there is a $t$ with positive exponent sign, denoted by $t_{l}$ for $1 \leq l \leq n$, which occurs both in a rule of kind $(1)$ and in a rule of kind $(\gamma)$, then in the construction of  $\Re''$ the rule of kind $(\gamma)$ has been replaced by a rule of kind $(5)$ in which there is no  $t$ with positive exponent sign. The same argument holds if  if there is a $t$ with positive exponent sign which occurs both in a rule of kind $(1)$ and in a rule of kind $(\delta)$, since $(\gamma)$ and $(\delta)$ have the same form (see proposition \ref{prop:R''_reduced}).
\end{proof}
\begin{lem}
\label{lem:t_moves}Let $w=S_{1}t_{1}S_{2}t_{2}..S_{k}t_{k}$ be a
non-empty word in $(X\cup\overline{X})^{*}$, with $S_{i}$ a sequence
of sources such that the length of $S_{i}$, $\ell(S_{i})$, satisfies
$0\leq\ell(S_{i})\leq\ell(w)$ for $1\leq i\leq k$. Then the application
of any rule in $\Re''$ on $w$ leads to one of the following three
cases for the $"t"$ which participated in the left-hand side of the
rule:

1. $t$ moves one step to the left and there is at most one source
added to $w$.

2. $t$ stays at the same place and there is one source added to $w$.

3. $t$ disappears and there is at most one source added to $w$.
In fact,
each application of a rule from class $(D)$ or of kind $t\overline{t} \rightarrow 1$
or $\overline{t}t \rightarrow 1$ on $w$ reduces the number
of $"t"$ by 2.
\end{lem}
\begin{proof}
The first case occurs when a rule from class $(A)$ or $(B)$ is applied,
the second case when a rule from class $(C)$ and the last case when
a rule from class $(D)$ or a rule of kind (0) is applied.
\end{proof}

 \begin{thm}
The rewriting system $\Re''$ is terminating.
\end{thm}

\begin{proof}
Let $w=S_{1}t_{1}S_{2}t_{2}..S_{k}t_{k}S_{k+1}$  be a word in $(X\cup\overline{X})^{*}$,
with  $0\leq n_{i}\leq \ell(w)$ for $1\leq i\leq k+1$.
We will show that after the application of any rule from  $\Re''$ on $w$, we obtain a word $w'$
 such that $V(w') <_{lex} V(w)$. We denote $w'=S'_{1}t'_{1}S'_{2}t'_{2}..S'_{k}t'_{k}S'_{k+1}$
where  $\ell(S'_{i})=n'_{i}$ for $1\leq i\leq k+1$.

If a rule from class (D) is applied on $w$, then
$V_{1}(w')=V_{1}(w)-2$ since each single application of a rule
from class $(D)$ on a word $w$ reduces the number of $"t"$ by 2.
This implies that  $V(w') <_{lex} V(w)$. The same holds for all
the rules of kind $(0)$: $t\overline{t} \rightarrow 1$
and $\overline{t}t \rightarrow 1$.\\

If a rule from class (C) is applied on $w$, then $V_{1}(w)=V_{1}(w')$ since the number of
occurrences of $t$ does not change but  $V_{2}(w')=V_{2}(w)-1$. This is due to the fact that
the $t^{+}$  on which the rule is applied is counted in $V_{2}(w)$
 (see remark \ref{rem:t+_neighX0}). So, $V(w') <_{lex} V(w)$.\\

From lemma \ref{lem:t_moves}, each application of a rule from class
$(A),(B)$ on $w$  keeps the number of $"t"$ the same, so $V_{1}(w)=V_{1}(w')$.
Moreover, the "$t^{+}$" which are counted in $V_{2}(w)$  cannot appear in the left-hand side,
 nor in the right-hand side of any rule from class $(A),(B)$, so  $V_{2}(w)=V_{2}(w')$.

First, assume a rule from class $(A)$ is applied on $w$ and let denote
by $t_{p}$  ($1\leq p \leq k$) the $t$ on which the rule is applied.
So, it holds that $n'_{i}=n_{i}$ for $i \neq p,p+1$,
$n'_{p}=n_{p}-1$ and $n'_{p+1}=n_{p+1}+2$.
The computation of $V_{3}(w')$ gives the following:\\
 $V_{3}(w')=\sum _{j=1}^{j=k}\sum_{i=0}^{i=k-j}2^{i}n'_{j}=
       \sum _{j=1}^{j=p-1}\sum_{i=0}^{i=k-j}2^{i}n'_{j} +
       \sum_{i=0}^{i=k-p}2^{i}n'_{p}+
       \sum_{i=0}^{i=k-p-1}2^{i}n'_{p+1}+
        \sum _{j=p+2}^{j=k}\sum_{i=0}^{i=k-j}2^{i}n'_{j}=
\sum_{j=1}^{j=p-1}\sum_{i=0}^{i=k-j}2^{i}n_{j} +
       \sum_{i=0}^{i=k-p}2^{i}(n_{p}-1)+
       \sum_{i=0}^{i=k-p-1}2^{i}(n_{p+1}+2)+
        \sum_{j=p+2}^{j=k}\sum_{i=0}^{i=k-j}2^{i}n_{j}=
          V_{3}(w) - 1(\sum_{i=0}^{i=k-p}2^{i}) +2(\sum_{i=0}^{i=k-p-1}2^{i})=
       V_{3}(w) - 2^{k-p} + 1(\sum_{i=0}^{i=k-p-1}2^{i})=
       V_{3}(w) - 2^{k-p} + (2^{k-p}-1)/(2-1)= V_{3}(w) - 1$\\
So,  $V_{3}(w')< V_{3}(w)$, which implies that $V(w')<_{lex} V(w)$. \\

Next, assume a rule from class $(B)$ is applied on $w$ and let denote
by $t_{p}$  ($1\leq p \leq k$) the $t$ on which the rule is applied.
So, it holds that $n'_{i}=n_{i}$ for $i \neq p,p+1$,
$n'_{p}=n_{p}-1$ and $n'_{p+1}=n_{p+1}+1$.
The computation of $V_{3}(w')$ gives the following:\\
 $V_{3}(w')=\sum_{j=1}^{j=k}\sum_{i=0}^{i=k-j}2^{i}n'_{j}=
          V_{3}(w) - 1(\sum_{i=0}^{i=k-p}2^{i}) +1(\sum_{i=0}^{i=k-p-1}2^{i})=
       V_{3}(w) - 2^{k-p} $\\
So,  $V_{3}(w')< V_{3}(w)$, which implies that $V(w')<_{lex} V(w)$. \\

At last, if a rule of kind $s\overline{s} \rightarrow 1$
or $\overline{s}s \rightarrow 1$ is applied on  $w$, then $V_{1}(w')=V_{1}(w)$,
$V_{2}(w')=V_{2}(w)$, since the number of occurrences of $t$ of any kind does not change.
In order to compare $V_{3}(w)$ and $V_{3}(w')$, we have to check several cases.\\

Assume there is at least one occurrence of $t$ in $w$.
If there is no occurrence of $t$ at the right of the subword on which
a rule of kind $s\overline{s} \rightarrow 1$
or $\overline{s}s \rightarrow 1$ is applied, i.e $s\overline{s}$ or $\overline{s}s$ is a subword of $S_{k+1}$, then $V_{3}(w')=V_{3}(w)$, since the computation of $V_{3}(w)$ does not take into account $n_{k+1}$, the length of $S_{k+1}$.
At last,  $V_{4}(w')=V_{4}(w)-2$, which implies $V(w')<_{lex} V(w)$. \\

Otherwise, if there is an occurrence of $t$, denoted by $t_{p}$, at the right of the subword
 on which a rule of kind $s\overline{s} \rightarrow 1$
or $\overline{s}s \rightarrow 1$ is applied, then $V_{3}(w')$ is not the same as $V_{3}(w)$, since $n'_{p}=n_{p}-2$ and $p \neq k+1$.
The computation of $V_{3}(w')$ gives the following:\\
 $V_{3}(w')=\sum _{j=1}^{j=k}\sum_{i=0}^{i=k-j}2^{i}n'_{j}=
          V_{3}(w) - 2(\sum_{i=0}^{i=k-p}2^{i}) $\\
So,  $V_{3}(w')< V_{3}(w)$, which implies that $V(w')<_{lex} V(w)$. \\

Assume there is no occurrence of $t$ in $w$, then $w=S_{1}$. Since the application of
a rule of kind $s\overline{s} \rightarrow 1$
or $\overline{s}s \rightarrow 1$ reduces the length of $S_{1}$ by 2, we have that  $n'_{1}=n_{1}-2$.
The computation of $V_{3}(w')$ gives the following:\\
 $V_{3}(w')=\sum _{j=1}^{j=1}\sum_{i=0}^{i=1-j}2^{i}n'_{j}=
          n'_{1}=n_{1}-2 = V_{3}(w)-2$\\
So, $V(w')<_{lex} V(w)$. \\

So, we have that after the application of any rule from  $\Re''$
on $w$, we obtain a word $w'$ which satisfies $V(w') <_{lex} V(w)$
and this implies that $\Re''$ is terminating .
\end{proof}

\subsection{The main result: Completeness of the rewriting system $\Re''$ }

\begin{thm}
The rewriting system $\Re''$ is finite and complete.
\end{thm}
\begin{proof}
$\Re''$ is terminating so it remains to show that $\Re''$ is locally
confluent, i.e that all the critical pairs, obtained from overlappings
between rules and from inclusion ambiguities of some rules in other
rules, resolve with no addition of new rules. This is done in the
lemmas in the next section, since the arguments used there are mostly
technical.
\end{proof}
The sets of rules in $\Re''$ are described in the lemmas
\ref{lem:add_rule_super_3} and \ref{lem:Add_rule_incl_over}. In
the following table we will give the list of critical pairs we
have to check, where ''-'' means that there is nothing to check
and $"\surd"$ means that there is an overlapping/inclusion which
is checked. A $"\surd"$ at row $i$ column $j$ means that there is
an overlapping of a rule of kind $i$ with a rule of kind $j$,
where the  rule of kind $i$ is at left and the  rule of kind $j$
is at right. As an example, inclusion ambiguities occur only
between rules of kind $(1)$ and rules of kind $(0)$. There can be
no overlapping between the rules of kind $(2')$, $(4)$, $(5)$,
$\alpha'$, $\beta$ and $\delta$, since in their left-hand side the
first letter is some "$s$" and the last letter is some "$t$".

\begin{tabular}{|c|c|c|c|c|c|c|c|c||c|}
\hline
&
0&
1&
2'&
3&
4&
5&
$\alpha$'&
$\beta$&
$\delta$\\
\hline
\hline
0&
-&
$\surd$&
$\surd$&
$\surd$&
$\surd$&
$\surd$&
$\surd$&
$\surd$&
$\surd$\\
\hline
1&
$\surd$&
-&
-&
-&
-&
-&
-&
-&
-\\
\hline
2'&
$\surd$&
-&
-&
$\surd$&
-&
-&
-&
-&
-\\
\hline
3&
$\surd$&
-&
-&
$\surd$&
-&
-&
-&
-&
-\\
\hline
4&
$\surd$&
-&
-&
$\surd$&
-&
-&
-&
-&
-\\
\hline
5&
$\surd$&
-&
-&
$\surd$&
-&
-&
-&
-&
-\\
\hline
$\alpha$'&
$\surd$&
-&
-&
-&
-&
-&
-&
-&
-\\
\hline
$\beta$&
$\surd$&
-&
-&
$\surd$&
-&
-&
-&
-&
-\\
\hline
\hline
$\delta$&
$\surd$&
-&
-&
-&
-&
-&
-&
-&
-\\
\hline
\end{tabular}

\subsection{Locally confluence of $\Re''$}

\begin{lem}
Critical pairs resulting from inclusion ambiguity of rule $(1)$ and
rule $(0)$ resolve.
\end{lem}
\begin{proof}
$\begin{array}{ccccccccc}
 &  & t_{l}\overline{t_{l}} &   &   and &  &   \overline{t_{l}}t_{l}\\
(1) & \swarrow &  \searrow & (0) &  & (1) & \swarrow   & \searrow \\
\overline{t_{j}}s_{k}\overline{t_{l}} &   & & 1 &  & \overline{t_{l}}\,\overline{t_{j}}s_{k} &  &   1\\
\downarrow & (2') &  &  &  & \downarrow & (3)\\
\overline{t_{j}}\,\overline{t_{l'}}s_{k'} &   &  &  & &\overline{s_{k'}}s_{k'} & \hookrightarrow\\
\downarrow & (3)\\
\overline{s_{k'}}s_{k'} & \hookrightarrow
\end{array}$
\end{proof}
\begin{lem}
Critical pairs resulting from overlapping between rule $(0)$ and
rule $(2')$ resolve.
\end{lem}
\begin{proof}
$\begin{array}{ccccccccc}
 &  & s_{k}\overline{t_{l}}t_{l} &  &   and &  &   \overline{s_{k}}s_{k}\overline{t_{l}}\\
(2') & \swarrow &   \searrow & (0) &  & (2') & \swarrow   & \searrow \\
\overline{t_{l'}}s_{k'}t_{l} &  &  &   s_{k} &   &\overline{s_{k}}\,\overline{t_{l'}}s_{k'} &  &
\overline{t_{l}}\\
\downarrow & (1)  &  &  & & \downarrow & (4)\\
\overline{t_{l'}}s_{k'}\overline{t_{j}}s_{k}   &  &  &  & & \overline{t_{l}}\,\overline{s_{k'}}s_{k'} &
\hookrightarrow\\
\downarrow & (2)\\
\overline{t_{l'}}t_{l'}s_{k} & \hookrightarrow
\end{array}$
\end{proof}
\begin{lem}
Critical pairs resulting from overlapping between rule $(0)$ and
rule $(3)$ resolve.
\end{lem}
\begin{proof}
$\begin{array}{ccccccccc}
 &  & t_{l}\overline{t_{l}}\,\overline{t_{j}} &  &   and &  &  \overline{t_{l}}\,\overline{t_{j}}t_{j}\\
(3) & \swarrow &   \searrow & (0) &  & (3) & \swarrow &   \searrow \\
t_{l}\overline{s_{k}} &  &  &   \overline{t_{j}} &  & \overline{s_{k}}\, t_{j}   &  &   \overline{t_{l}}\\
\downarrow & (1) &  &  &    & \downarrow & (1)\\
\overline{t_{j}}s_{k'}\overline{s_{k}} &  &  &  &  & \overline{s_{k}}\,\overline{t_{l'}}s_{k'}\\
\downarrow & (0) &  &  &   &\downarrow & (4)\\
\overline{t_{j}} & \hookrightarrow &  &  &  & \overline{t_{l}}\,\overline{s_{k'}}s_{k'} & \hookrightarrow\end{array}$
\end{proof}
\begin{lem}
Critical pairs resulting from overlapping between rule $(0)$ and
rule $(4)$ resolve.
\end{lem}
\begin{proof}
$\begin{array}{ccccccccc}
 &  & s_{k}\overline{s_{k}}\,\overline{t_{l'}} &  &   and &  &   \overline{s_{k}}\,\overline{t_{l'}}t_{l'}\\
(4) & \swarrow &  \searrow & (0) &  & (4) & \swarrow &   \searrow \\
s_{k}\overline{t_{l}}\,\overline{s_{k'}} &  &  &   \overline{t_{l'}}{} &  & \overline{t_{l}}\,\overline{s_{k'}}\,
t_{l'} &  &     \overline{s_{k}}\\
\downarrow & (2') &  &  &  &   \downarrow & (1)\\
\overline{t_{l'}}s_{k'}\overline{s_{k'}} & \hookrightarrow &  &  &  &
\overline{t_{l}}\,\overline{s_{k'}}\,\overline{t_{q}}s_{q'}\\
 & (0) &  &  &  & \downarrow & (4)\\
 &  &  &  &    &\overline{t_{l}}\,\overline{t_{j}}\,\overline{s_{q'}}s_{q'} & \hookrightarrow & (3)\end{array}$
\end{proof}
\begin{lem}
Critical pairs resulting from overlapping between rule $(0)$ and
rule $(5)$ resolve.
\end{lem}
\begin{proof}
$\begin{array}{ccccccccc}
 &  & s_{m}\overline{s_{m}}\overline{t_{j}} &  &   and &  &   \overline{s_{m}}\overline{t_{j}}t_{j}\\
(5) & \swarrow &   \searrow & (0) &  & (5) & \swarrow &   \searrow \\
s_{m}\overline{t_{p}}s_{m'}\overline{s_{k}}&  &  &  \overline{ t_{j}} & &  \overline{t_{p}}\,
s_{m'}\overline{s_{k}}t_{j} &  &    \overline{s_{m}}\\
\downarrow & (\alpha') &  &  &  &   \downarrow & (1)\\
\overline{t_{j}}s_{k}\overline{s_{m'}}s_{m'}\overline{s_{k}} & \hookrightarrow &  &  &  &
\overline{t_{p}}s_{m'}\overline{s_{k}}\,\overline{t_{l'}}s_{k'}\\
 &  &  &  &  &   \downarrow & (4)\\
 &  &  &  &  &   \overline{t_{p}}s_{m'}\overline{t_{l}}\, \overline{s_{k'}}s_{k'}\\
 &  &&  &  &  \downarrow & (0),(3)\\
&  &  &  &  &\overline{t_{p}}t_{p}\overline{s_{m}}\hookrightarrow
\end{array}$
\end{proof}
\begin{lem}
Critical pairs resulting from overlapping between rule  $(0)$ and
rule $(\alpha')$ resolve.
\end{lem}
\begin{proof}
$\begin{array}{ccccccccc}
 &  & s_{m}\overline{t_{p}}\, t_{p} &  &   and &  &  \overline{s_{m}}s_{m}\overline{t_{p}}\\
(\alpha') & \swarrow &   \searrow & (0) &  & (\alpha') & \swarrow &   \searrow \\
\overline{t_{j}}s_{k}\overline{s_{m'}}t_{p} &  &  &   s_{m} &  &
\overline{s_{m}}\,\overline{t_{j}}s_{k}\overline{s_{m'}} &  &     \overline{t_{p}}\\
\downarrow & (\delta) &  &  &  &   \downarrow & (5)\\
\overline{t_{j}}s_{k}\overline{t_{l}}s_{m} &  &  &  & &
\overline{t_{p}}s_{m'}\overline{s_{k}}s_{k}\overline{s_{m'}} &   \hookrightarrow\\
\downarrow & (2')\\
\overline{t_{j}}\,\overline{t_{l'}}s_{k'}s_{m} \\
\downarrow & (3)\\
 \overline{s_{k'}}s_{k'}s_{m} & \hookrightarrow
\end{array}$
\end{proof}
\begin{lem}
Critical pairs resulting from overlapping between rule $(0)$
and rule  $(\beta)$  resolve.
\end{lem}
\begin{proof}
$\begin{array}{ccccccccc}
 &  & s_{m'}\overline{t_{l}}\, t_{l} &  &   and   &  & \overline{s_{m'}}s_{m'}\overline{t_{l}}\\
(\beta) & \swarrow &   \searrow & (0) &  & (\beta) & \swarrow &   \searrow & (0)\\
t_{p}\overline{s_{m}}t_{l} &  &  &   s_{m'} &  & \overline{s_{m'}}t_{p}\overline{s_{m}} &  &
\overline{t_{l}}\\
\downarrow & (1) &  &  &  &   \downarrow & (\delta)\\
t_{p}\overline{s_{m}}\,\overline{t_{j}}s_{k} &  &  &  &  &  \overline{t_{l}}s_{m}\overline{s_{m}} &  &
\hookrightarrow\\
\downarrow & (5)\\
\overline{t_{p}}\,\overline{t_{p}}s_{m'}\overline{s_{k}}s_{k} &  &   \hookrightarrow\end{array}$
\end{proof}
\begin{lem}
Critical pairs resulting from overlapping between rule $(0)$
and rule  $(\delta)$  resolve.
\end{lem}
\begin{proof}
$\begin{array}{ccccccccc}
 &  & s_{m'}\overline{s_{m'}}t_{p} &  &   and &  &   \overline{s_{m'}}t_{p}\,\overline{t_{p}}\\
(\delta) & \swarrow &   \searrow & (0) &  & (\delta) & \swarrow &   \searrow \\
s_{m'}\overline{t_{l}}s_{m} &  &  &   t_{p} &  & \overline{t_{l}}\, s_{m}\overline{t_{p}} &  &
\overline{s_{m'}}\\
\downarrow & (\beta) &  &  &  &  \downarrow & (\alpha')\\
t_{p}\overline{s_{m}}s_{m} & \hookrightarrow &  &  &  &
\overline{t_{l}}\,\overline{t_{j}}s_{k}\overline{s_{m'}}\\
   &  &  &  &  & \downarrow & (3)\\
   &  &  &  &  & \overline{s_{k}}\, s_{k}\overline{s_{m'}} &  & \hookrightarrow
\end{array}$
\end{proof}

\begin{lem}
Critical pairs resulting from overlapping between rule $(2')$ and
rule $(3)$ resolve.
\end{lem}
\begin{proof}
$\begin{array}{cccccccc}
 &  & s_{k}\overline{t_{l}}\,\overline{t_{j}}\\
(2') & \swarrow &  & \searrow & (3)\\
\overline{t_{l'}}s_{k'}\overline{t_{j}} &  &  &  & s_{k}\overline{s_{k}}{}\\
\downarrow & (2) &  &  & \downarrow\\
\overline{t_{l'}}t_{l'} & \hookrightarrow &  &  & 1\\
\\\end{array}$
\end{proof}
\begin{lem}
Critical pairs resulting from overlapping between rule $(3)$ and
rule $(3)$ resolve.
\end{lem}
\begin{proof}
$\begin{array}{cccccccc}
 &  & \overline{t_{l}}\,\overline{t_{j}}\,\overline{t_{l'}}\\
(3) & \swarrow &  & \searrow & (3)\\
\overline{s_{k}}\,\overline{t_{l'}} & \hookrightarrow &  &  & \overline{t_{l}}\,\overline{s_{k'}}{}\\
 & (4) &  &  & \\
\\\end{array}$
\end{proof}
\begin{lem}
Critical pairs resulting from overlapping between rule $(5)$ and
rule $(3)$ resolve.
\end{lem}
\begin{proof}
$\begin{array}{cccccccc}
 &  & \overline{s_{m}}\,\overline{t_{j}}\,\overline{t_{l'}}\\
(5) & \swarrow &  & \searrow & (3)\\
\overline{t_{p}}s_{m'}\overline{s_{k}}\,\overline{t_{l'}} &  &  &  & \overline{s_{m}}\,\overline{s_{k'}}\\
\downarrow & (4)\\
\overline{t_{p}}s_{m'}\overline{t_{l}}\,\overline{s_{k'}} &  &  &  & 1\\
\downarrow & (\beta)\\
\overline{t_{p}}t_{p}\overline{s_{m}}\,\overline{s_{k'}} &  & \hookrightarrow\end{array}$
\end{proof}
\begin{lem}
Critical pairs resulting from overlapping between rule $(\beta)$
and rule $(3)$ resolve.
\end{lem}
\begin{proof}
$\begin{array}{cccccccc}
 &  & s_{m'}\overline{t_{l}}\,\overline{t_{j}}\\
(\beta) & \swarrow &  & \searrow & (3)\\
t_{p}\overline{s_{m}}\,\overline{t_{j}} &  &  &  & s_{m'}\,\overline{s_{k}}\\
\downarrow & (5)\\
t_{p}\overline{t_{p}}s_{m'}\overline{s_{k}} &  & \hookrightarrow\\
\\\end{array}$
\end{proof}
Using the same argument as in claim \ref{cla:if_setofm_then_setof m'},
we obtain that if $\Phi_{m'}$ is in $\Re'$ then the set of rules
$\Phi_{m''}$ is also in $\Re'$, where $\Phi_{m''}$ is:
$\begin{array}{c}
(1_{l',q,q'}):t_{l'}\rightarrow\overline{t_{q}}s_{q'} \\
(2_{l',q,q'}):s_{q'}\overline{t_{l'}}\rightarrow t_{q} \\
(3_{l',q,q'}):\overline{t_{l'}}\,\overline{t_{q}}\rightarrow\overline{s_{q'}} \\
\end{array}$ and the following rule of kind $(4):\overline{s_{k'}}\,\overline{t_{q}}\rightarrow\overline{t_{j}}\,
\overline{s_{q'}}$ is added to $\Re'$.

\begin{lem}
Critical pairs resulting from overlapping between rule $(4)$ and
rule $(3)$ resolve.
\end{lem}
\begin{proof}
$\begin{array}{cccccccc}
 &  & \overline{s_{k}}\,\overline{t_{l'}}\,\overline{t_{q}}\\
(4) & \swarrow &  & \searrow & (3)\\
\overline{t_{l}}\,\overline{s_{k'}}\,\overline{t_{q}} &  &  &  & \overline{s_{k}}\,\overline{s_{q'}}{}\\
\downarrow & (4)\\
\overline{t_{l}}\,\overline{t_{j}}\,\overline{s_{q'}} &  & \hookrightarrow\\
 &  & (3)\\
\end{array}$
\end{proof}


\begin{thebibliography}{10}
\bibitem[1]{1}R.V. Book and F. Otto, String-Rewriting systems, Springer-Verlag,
1993.
\bibitem[2]{2}F.Chouraqui, Rewriting Systems in Alternating Knot Groups, Int. Journal of Alg.and Comput.,vol.16,Aug2006.
\bibitem[3]{3}R.H. Crowell and R.H. Fox, Introduction to Knot Theory, Blaisdell
Pub.Co., 1965.
\bibitem[4]{4}D.E. Cohen, String Rewriting-A Survey for Group Theorists, in Geometric
Group Theory, vol.1, G.A.Niblo and M.A.Roller (eds), Lecture Notes
of the London .Math.Soc.181 (1993).
\bibitem[5]{5}S. Hermiller, J. Meier, Artin groups, rewriting systems and three-manifolds,
J.Pure Appl.Algebra 136 (1999), 141-156.
\bibitem[6]{6}P. Le Chenadec, Canonical forms in finitely presented algebras, Pitman,
London, 1986.
\bibitem[7]{7}R.C. Lyndon and P.E. Schupp, Combinatorial Group Theory, Springer-Verlag,
1977.
\bibitem[8]{8}C.C.Squier, Word prolems and an homological finiteness condition for monoids,
J.Pure Appl.Algebra 49 (1987), pp.201-217.
\bibitem[9]{9}C.M.Weinbaum, The word and conjugacy problems for the knot group of
any tame,prime,alternating knot, Proc.Amer.Math.Soc.30 (1971), pp.22-26.
\bibitem[10]{10}R.J.Wilson, Introduction to Graph Theory, Oliver and Boyd, Edinburgh,1972.

\end{thebibliography}
\end{document}